\documentclass[conference]{IEEEtran}
\IEEEoverridecommandlockouts
\usepackage{tikz}
\usetikzlibrary{arrows.meta, positioning}
\usepackage{cite}
\usepackage{amsmath,amssymb,amsfonts}
\usepackage{algorithmic}
\usepackage{graphicx}
\usepackage{textcomp}
\usepackage{xcolor}
\usepackage{hyperref}
\usepackage{url}
\usepackage{bm} 
\def\BibTeX{{\rm B\kern-.05em{\sc i\kern-.025em b}\kern-.08em
    T\kern-.1667em\lower.7ex\hbox{E}\kern-.125emX}}
\begin{document}

\title{Influence Diagrams for\\
Robust Multi-Target Tracking\\

}

\author{\IEEEauthorblockN{ Priyank Behera}
\IEEEauthorblockA{\textit{Department Of Computer Science} \\
\textit{Purdue University}\\
West Lafayette, US \\
beherap@purdue.edu}
\and
\IEEEauthorblockN{C. Robert Kenley}
\IEEEauthorblockA{\textit{Edwardson School of Industrial Engineering} \\
\textit{Purdue University}\\
West Lafayette, US \\
kenley@purdue.edu}
\and
}

\maketitle

\begin{abstract}
Multi-Target Tracking (MTT) is foundational for radar, defense, and autonomous systems, where tracking accuracy directly affects decision-making and safety. For linear systems with Gaussian process and measurement noise, the Kalman filter remains the gold standard for state estimation. However, its performance can degrade in real-world scenarios where measurement noise is temporally correlated. This violates the Kalman filter's white-noise assumptions. Various approaches include state augmentation of the Kalman filter, but this approach is susceptible to failure due to ill-conditioned problem formulations. This work investigates the limitations of classical Kalman filtering in colored noise environments and presents an influence diagram-based approach to the Joint Probabilistic Data Association Filter (JPDAF). Simulation results on benchmark scenarios demonstrate that the Influence Diagram JPDAF (ID-JPDAF) achieves lower root mean square error (RMSE) than classical methods. These findings highlight the potential of influence diagram models for advancing multi-target tracking performance in radar and related applications.
\end{abstract}


\section{Introduction}
Multi-target tracking (MTT) is a core problem in radar, surveillance, and autonomous systems. These systems require reliable estimation of multiple moving objects in noisy and cluttered environments. For the single target case, the Kalman filter \cite{kalman1960} and its variants remain widely adopted for such tasks due to their computational efficiency and optimality under linear-Gaussian assumptions. However, classical
Kalman filters exhibit two critical limitations that can significantly degrade estimation accuracy or, in adverse conditions, result in filter divergence and tracking failure. First, the Kalman filter assumes that measurement noise is independently and identically distributed (i.i.d.) Gaussian random variables. In practice, the i.i.d. assumption can fail because of the temporal correlation between noise components, also known as colored noise. Second, the Kalman filter involves the inversion of a matrix, specifically \textit{innovation covariance}. This can cause the system to become ill-conditioned, resulting in a numerically unstable inversion, high root mean squared error (RMSE), and poor filter performance \cite{maybeck1982stochastic}. 

Within the context of MTT, there are various algorithms in the radar world. Some examples include the Nearest-Neighbor Filter \cite{blackman1986multiple}, the Joint Probabilistic Data Association Filter (JPDAF) \cite{chang1984joint},
Multiple Hypothesis Tracking \cite{reid2003algorithm}, and the Interacting Multiple-Model Filter \cite{blackman1986multiple}. These methods use the Kalman filter, or one of its variants, to estimate the state of all the objects in question. 

This work focuses on JPDAF due to its soft measurement update mechanism, which balances accuracy and complexity in cluttered environments. We systematically evaluate the performance of the classical Kalman filter and the influence diagram-based Kalman filter~\cite{kenley1986influence} within the JPDAF framework for multi-target tracking. Through simulation experiments, we characterize the estimation accuracy and robustness of both approaches, highlighting the degradation of standard filtering methods under colored noise and filter mismatch. The results demonstrate that the influence diagram-based filter achieves superior performance in challenging scenarios, offering improved numerical stability and resilience to model assumptions.

\vspace{0.5em}
\section{Background}
\subsection{State-Space Model}
MTT extends the classical state-space formulation to simultaneously estimate the trajectories of $N$ independent or interacting targets in cluttered environments. The joint state-space model for $N$ targets at time $k$ can be represented as \cite{bryson_chapter_1975}:
\begin{align}
    \mathbf{X}_{k} &= \mathbf{F}_{k-1} \mathbf{X}_{k-1} + \mathbf{W}_{k} \label{eq:mtt-state-evolution} \\
    \mathbf{Z}_{k} &= \mathbf{H}_{k} \mathbf{X}_{k} + \mathbf{V}_{k} \label{eq:mtt-measurement}
\end{align}
where:
\begin{itemize}
    \item $\mathbf{X}_{k} = [\mathbf{x}_{k}^{(1)\top},\, \ldots,\, \mathbf{x}_{k}^{(N)\top}]^\top \in \mathbb{R}^{N n}$ is the stacked state vector of all $N$ targets, with each target state $\mathbf{x}_k^{(i)} \in \mathbb{R}^n$.
    \item $\mathbf{Z}_{k} = [\mathbf{z}_{k}^{(1)\top},\, \ldots,\, \mathbf{z}_{k}^{(N)\top}]^\top \in \mathbb{R}^{N m}$ denotes the stacked measurement vector, where $\mathbf{z}_k^{(j)} \in \mathbb{R}^m$ is the measurement associated with the $j$-th target at time $k$.
    \item $\mathbf{F}_{k-1} \in \mathbb{R}^{N n \times N n}$ is the diagonal transition matrix, with each block corresponding to individual target dynamics.
    \item $\mathbf{H}_k \in \mathbb{R}^{N m \times N n}$ is the block-diagonal measurement matrix, mapping the joint state to the measurement space.
    \item $\mathbf{W}_{k} \sim \mathcal{N}(\mathbf{0}, \mathbf{Q}_k)$ is the joint process noise, where $\mathbf{Q}_k$ is block-diagonal (if targets are independent).
    \item $\mathbf{V}_k \sim \mathcal{N}(\mathbf{0}, \mathbf{R}_k)$ is the joint diagonal measurement noise
\end{itemize}

This multi-target formulation provides a unified framework for representing the evolution and observation of multiple targets. In practice, data association must be resolved to correctly pair measurements to targets, and probabilistic association filters (such as JPDAF) are commonly employed to address measurement origin uncertainty and false alarms.

\subsection{Joint Probabilistic Data Association Filter (JPDAF)}
In MTT, associating measurements to the correct targets amidst clutter and missed detections is a central challenge. JPDAF addresses this by updating each track with a Bayesian mixture over all feasible measurement-to-track associations~\cite{barshalom1993}.

The track prediction is:
\begin{align}
    {\mathbf{x}}_{k|k-1}^{(i)} &= \mathbf{F}_{k-1} {\mathbf{x}}_{k-1|k-1}^{(i)} \\
    \mathbf{P}_{k|k-1}^{(i)} &= \mathbf{F}_{k-1} \mathbf{P}_{k-1|k-1}^{(i)} \mathbf{F}_{k-1}^\top + \mathbf{Q}_{k-1}
\end{align}

At time $k$, let $\{\mathbf{z}_k^{(j)}\}_{j=1}^{M_k}$ be the set of $M_k$ measurements and consider $N$ tracks. For track $i$, define association probabilities $\beta_k^{(i,j)}$ for each measurement $j$ and $\beta_k^{(i,0)}$ for a missed detection, so that
   $ \sum_{j=0}^{M_k} \beta_k^{(i,j)} = 1.$
These are computed via Bayes' rule:
\begin{align}
    \beta_k^{(i,j)} = \frac{P(\mathbf{z}_k^{(j)}|{\mathbf{x}}_{k|k-1}^{(i)})P(\mathcal{A}_k^{(i,j)})}
    {\sum_{l=0}^{M_k} P(\mathbf{z}_k^{(l)}|{\mathbf{x}}_{k|k-1}^{(i)})P(\mathcal{A}_k^{(i,l)})},
\end{align}
where $P(\mathcal{A}_k^{(i,j)})$ is the prior (e.g., $P_D P_G / \lambda$ for detection, $1-P_D P_G$ for missed), and $P(\mathbf{z}_k^{(j)}|{\mathbf{x}}_{k|k-1}^{(i)})$ is the Kalman likelihood,
$    P(\mathbf{z}_k^{(j)}|{\mathbf{x}}_{k|k-1}^{(i)}) = 
    \mathcal{N}(\mathbf{z}_k^{(j)};\, \mathbf{H}_k {\mathbf{x}}_{k|k-1}^{(i)},\, \mathbf{S}_k^{(i)}),
$
Switch innovation covariance $\mathbf{S}_k^{(i)} = \mathbf{H}_k \mathbf{P}_{k|k-1}^{(i)} \mathbf{H}_k^\top + \mathbf{R}_k$.

The track update is a weighted sum:
\begin{align}
    {\mathbf{x}}_{k|k}^{(i)} &= {\mathbf{x}}_{k|k-1}^{(i)} + \mathbf{K}_k^{(i)} \left(\sum_{j=1}^{M_k} \beta_k^{(i,j)} \mathbf{y}_k^{(i,j)}\right),\\
    \mathbf{P}_{k|k}^{(i)} &= \beta_k^{(i,0)} \mathbf{P}_{k|k-1}^{(i)} + (1 - \beta_k^{(i,0)}) \tilde{\mathbf{P}}_{k|k}^{(i)} \notag \\
    &\quad + \sum_{j=1}^{M_k} \beta_k^{(i,j)} \mathbf{e}_k^{(i,j)}\mathbf{e}_k^{(i,j)\top} - \mathbf{e}_k^{(i)}\mathbf{e}_k^{(i)\top},
\end{align}
where $\mathbf{y}_k^{(i,j)} = \mathbf{z}_k^{(j)} - \mathbf{H}_k {\mathbf{x}}_{k|k-1}^{(i)}$, $\mathbf{e}_k^{(i,j)} = \mathbf{K}_k^{(i)} \mathbf{y}_k^{(i,j)}$, and $\mathbf{e}_k^{(i)} = \mathbf{K}_k^{(i)} \sum_{j=1}^{M_k} \beta_k^{(i,j)} \mathbf{y}_k^{(i,j)}$, and $\tilde{\mathbf{P}}^{(i)}_{k|k} = (\mathbf{I} - \mathbf{K}^{(i)}_k \mathbf{H}_k)\, \mathbf{P}^{(i)}_{k|k-1}$. 
This Bayesian updating allows JPDAF to robustly handle ambiguous data associations in moderately cluttered environments at a good computational cost.

\vspace{1em}
\subsection{Colored Noise and State Augmentation}
Classical Kalman filters assume that the process and measurement noise are white, i.e., temporally uncorrelated. However, in many radar and sensing applications measurement noise exhibit colored noise. Colored noise often arises from sensor memory, environmental effects, or system imperfections. If unaddressed, colored noise can lead to biased estimates and degraded tracking performance. A standard approach for handling colored measurement noise is to \emph{augment} the state-space model so that the noise becomes part of the extended system state. Let the measurement noise $\mathbf{v}_k$ follow:
\begin{align}
    \mathbf{v}_k = \rho \mathbf{v}_{k-1} + \bm{\xi}_k,
\end{align}
where $|\rho| < 1$ is the coloredness factor, and the zero-mean $\mathbf{\xi}_k \sim \mathcal{N}(\mathbf{0}, \mathbf{\sigma^2I})$ is white Gaussian noise. The measurement equation is then:
\begin{align}
    \mathbf{z}_k = \mathbf{H}_k \mathbf{x}_k + \mathbf{v}_k
\end{align}

By defining an augmented state vector
$\tilde{\mathbf{x}}_k = [\mathbf{x}_k^\top,\, \mathbf{v}_k^\top]^\top$

$\in \mathbb{R}^{n + m}$
the system dynamics and measurement equations become
\begin{align}
    \tilde{\mathbf{x}}_k &= 
    \underbrace{
        \begin{bmatrix}
            \mathbf{F}_{k-1} & \mathbf{0} \\
            \mathbf{0} & \rho \mathbf{I}
        \end{bmatrix}
    }_{\tilde{\mathbf{F}}_{k-1}}
    \tilde{\mathbf{x}}_{k-1}
    +
    \begin{bmatrix}
        \mathbf{w}_k \\
        \bm{\xi}_k
    \end{bmatrix} \\
    \mathbf{z}_k &= 
    \underbrace{
        \begin{bmatrix}
            \mathbf{H}_k & \mathbf{I}
        \end{bmatrix}
    }_{\tilde{\mathbf{H}}_k}
    \tilde{\mathbf{x}}_k
\end{align}

The process noise for the augmented model is
$
\tilde{\mathbf{w}}_k = \begin{bmatrix}
    \mathbf{w}_k \\
    \bm{\xi}_k
\end{bmatrix},
$

This augmented formulation allows standard Kalman prediction and update equations to be applied to the augmented state $\tilde{\mathbf{x}}_k$ and its associated covariance.
\vspace{2em}
\section{Influence Diagram-Based JPDAF (ID-JPDAF)}
\vspace{1em}
The classical JPDAF utilizes a Kalman filter for each of the targets that are being tracked represented by a stacked vector of \textit{N} targets. Then, it performs a joint probabilistic data association across all measurement-to-track associations. Finally, the measurements are weighted with their association probabilities to achieve the final predicted state. 

\begin{figure}[ht]
    \centering
    \begin{tikzpicture}[
        node distance=0.8cm and 1.05cm,
        block/.style={draw, rounded corners, align=center, font=\scriptsize, minimum width=2.1cm, minimum height=0.5cm, fill=gray!10},
        datain/.style={draw, rounded corners, align=center, font=\scriptsize, minimum width=1.6cm, minimum height=0.45cm, fill=blue!12},
        arrow/.style={-Stealth, thick, shorten >=1pt, shorten <=1pt}
    ]
    \node[datain] (prev) {Prior State};
    \node[block, right=of prev] (kf) {Kalman\\Prediction};
    \node[block, below=of kf] (assoc) {Data Association\\(Compute $\beta$)};
    \node[block, below=of assoc] (update) {Kalman\\Update};
    \node[block, below=of update] (weight) {Posterior\\State};
    \node[datain, right=of assoc] (meas) {New\\Measurements};

    \draw[arrow] (prev) -- (kf);
    \draw[arrow] (kf) -- (assoc);
    \draw[arrow] (assoc) -- (update);
    \draw[arrow] (update) -- (weight);
    \draw[arrow] (meas) -- (assoc);

    \draw[arrow, dashed] (weight.south) -- ++(0,-0.28) -| (prev.south);

    \end{tikzpicture}
    \caption{JPDAF Pipeline}
    \label{fig:jpda_block_diagram}
\end{figure}
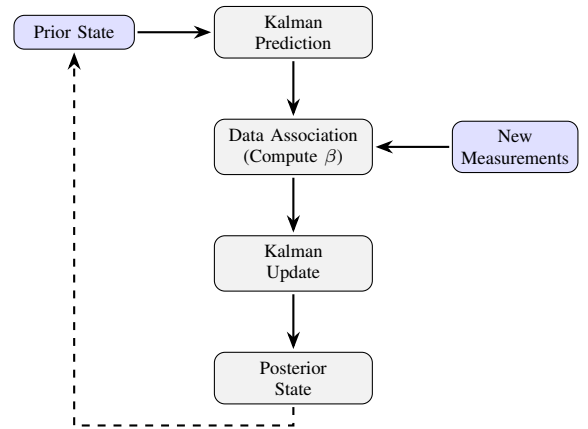

Fig.~\ref{fig:jpda_block_diagram} illustrates the general architecture of the JPDAF pipeline. The primary methodological contributions of this work are focused on the modification of the Kalman prediction and update components, in which the influence diagram framework is incorporated to inform the filtering process. In addition, several minor adjustments are made to the data association stage, including an updated measurement gating strategy and an alternative approach to the computation of the Mahalanobis distance. The filtering algorithms implemented in this work are derived from the discrete-time influence diagram formulation described in~\cite{kenley1986influence}. Kovacich successfully implemented these algorithms for multitarget tracking \cite{michael_kovacich_application_1990}. 

\subsection{Influence Diagrams}

Influence diagrams, also known as Bayesian networks, are probabilistic graphical models that provide a framework for inference~\cite{howard2005influence}. Traditionally, most practical influence diagram implementations were restricted to discrete variables. However, Kenley~\cite{kenley1986influence} extended this framework to continuous variables, introducing the \emph{Gaussian} influence diagram. In this framework, each node represents a Gaussian random variable, while directed arcs encode conditional linear dependencies between variables and the propagation of information.

\begin{figure}[t]
    \centering
    \begin{tikzpicture}[%
        node distance=0.95cm and 1.5cm,
        >=Stealth,
        every node/.style={
            align=center,
            text width=1.35cm,
            font=\normalsize
        },
        bigc/.style={circle, draw=black, thick, fill=red!18, minimum size=1.7cm},
        bigo/.style={circle, draw=black, thick, fill=cyan!18, minimum size=1.7cm}
    ]
    \node[bigc] (x1)  {$\mathit{x}_1$\\[0.1em] {\footnotesize $({x}_{1|0},\,P_{1|0})$}};
    \node[bigo, below=of x1] (z1) {$\mathit{z}_1$\\[0.1em] {\footnotesize $(H{x}_{1|0},\,0)$}};
    \node[bigo, right=of x1] (x2) {$\mathit{x}_2$\\[0.1em] {\footnotesize $({x}_{2|0},\,0)$}};
    \node[bigc, above=of x2] (w2) {$\mathit{w}_2$\\[0.1em] {\footnotesize $(0,\,Q_2)$}};
    \node[bigc, left=of z1] (v1) {$\mathit{v}_1$\\[0.1em] {\footnotesize $(0,\,R_1)$}};
    \draw[->, thick] (x1) -- node[above, yshift=1pt] {\footnotesize$F_1$} (x2);
    \draw[->, thick] (x1) -- node[left, xshift=10pt, yshift=0pt] {\footnotesize$H_1$} (z1);
    \draw[->, thick] (v1) -- (z1);
    \draw[->, thick] (w2) -- (x2);
    \end{tikzpicture}
    \caption{Gaussian influence diagram for discrete-time filtering. Each of these nodes is assumed to be normally distributed, taking the form $x \sim \mathcal{N}(\mu,\, \Sigma)$. The blue nodes are deterministic, and the pink nodes are stochastic
}
    \label{fig:gaussian_id}
\end{figure}
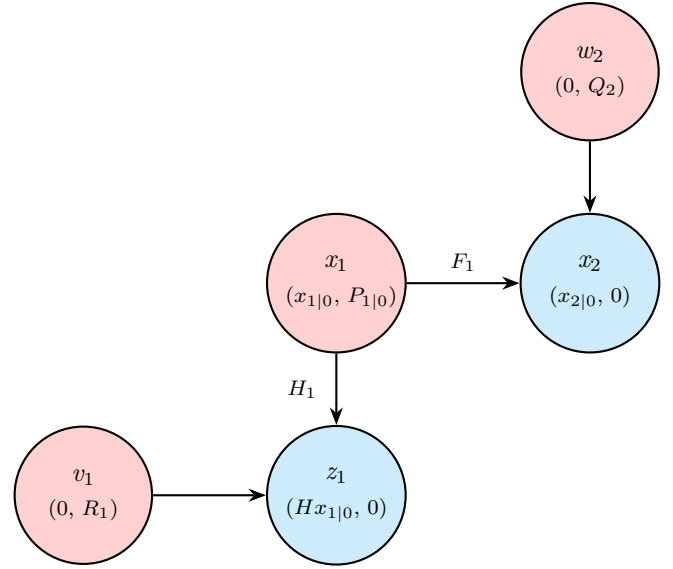

As shown in Fig.~\ref{fig:gaussian_id}, the Gaussian influence diagram facilitates both probabilistic inference and enables intuitive graphical manipulations such as node removal and arc reversal, corresponding to marginalization and conditioning in Bayesian inference. Instead of a covariance matrix $P$, the influence diagram performs operations on $V$ and $B$~\cite{kenley1986influence}, \cite{shachter_gaussian_1989}. Each node $X_j$ is parameterized by a conditional mean and variance given its parents $X_{C(j)}$. The influence of parent $X_k \in C(j)$ on $X_j$ is encoded by the regression coefficient $b_{kj}$, and the conditional variance (or covariance matrix) is denoted by $v_j$.

A multivariate normal random vector $\mathbf{x} \sim \mathcal{N}(\boldsymbol{\mu}, \mathbf{X})$ can be equivalently represented in Gaussian influence diagram form by decomposing the covariance matrix $\mathbf{X}$ into (i) a strictly upper triangular matrix $\mathbf{B}$ of arc (regression) coefficients and (ii) a vector $\mathbf{V}$ of conditional variances:
\begin{align}
    x_j &= \sum_{k=1}^{j-1} B_{k j} x_k + \epsilon_j, \quad \epsilon_j \sim \mathcal{N}(0, v_j), \quad j=1, \ldots, n \\
    \mathbf{x} &\sim \mathcal{N}(\boldsymbol{\mu}, \mathbf{X}) \iff (\boldsymbol{\mu}, \mathbf{B}, \mathbf{V})
\end{align}
Here, $B_{k j}$ and $V_j$ are computed recursively from the entries of $\mathbf{X}$:
\begin{align}
    B_{k j} &= \left[\mathbf{P}_{1:(j-1),1:(j-1)} \, \mathbf{X}_{1:(j-1),j}\right]_k \\
    V_j &= X_{j j} - \mathbf{X}_{j,1:(j-1)} \mathbf{B}_{1:(j-1),j}
\end{align}
where $\mathbf{P}$ is the inverse (or generalized inverse) of $\mathbf{X}$. This conversion allows us to perform operations such as arc reversal and arc removal, using $\mathbf{B}$ and $\mathbf{V}$. 

\subsection{Kalman Prediction}

The time update for the Gaussian influence diagram propagates the state from time $k-1$ to $k$. This is achieved by augmenting the current state and process noise in block form, converting the process noise covariance to influence diagram parameters, and then marginalizing the intermediate variables via node removal~\cite{kenley1986influence}. $\mathbf{B_{k-1|k-1}}$ and $\mathbf{V_{k-1|k-1}}$ are derived from the prior covariance $\mathbf{P_{k-1|k-1}}$. The key steps are:
\begin{enumerate}
    \item \textbf{Process noise conversion:} Convert the process noise covariance $\mathbf{Q_k}$ into influence diagram parameters $(\mathbf{B_q}, \mathbf{V_q})$:
    \begin{equation}
        \mathbf{Q_k} \longrightarrow (\mathbf{B_q},\mathbf{V_q})
    \end{equation}
    \item \textbf{Augmented influence diagram:} Form an augmented block system combining state and process noise:
    \begin{equation}
        \mathbf{B_{\text{aug}}} = 
        \begin{bmatrix}
            \mathbf{B_q} & \mathbf{0} & \mathbf{I} \\
            \mathbf{0} & \mathbf{B_{k-1|k-1}} & 
            \mathbf{F^\top} \\
            \mathbf{0} & \mathbf{0} & \mathbf{0}
        \end{bmatrix}, \quad
        \mathbf{V_{\text{aug}}} = \begin{bmatrix} \mathbf{V_q} \\ \mathbf{V_{k-1|k-1}} \\ \mathbf{0} \end{bmatrix}
    \end{equation}
    \item \textbf{Node removal:} Marginalize the process noise and prior state nodes to yield the updated parameters for the predicted state:
    \begin{equation}
        (\mathbf{B_{k|k-1}}, \mathbf{V_{k|k-1}}) = \texttt{removal}(\mathbf{B_{\text{aug}}}, \mathbf{V_{\text{aug}}})
    \end{equation}
    \item \textbf{Mean update:} The predicted mean is updated by the transition model:
    \begin{equation}
        \mathbf{u_{k|k-1}} = \mathbf{F u_{k-1|k-1}}
    \end{equation}
\end{enumerate}

Mathematically, node removal updates the remaining parameters as follows:
\begin{align}
    \mathbf{V}_{\text{rem}} &= \mathbf{V} + \mathbf{B}_{*j}^2\, \mathbf{V}_j, \\
    \mathbf{B}_{\text{rem}} &= \mathbf{B}_{*} + \mathbf{B}_{*j} \mathbf{B}_{j*},
\end{align}
where the subscript \( j \) indicates the node being removed and \( * \) indexes the remaining nodes.

This procedure produces the predicted influence diagram parameters \( (\mathbf{u}_{k|k-1}, \mathbf{B}_{k|k-1}, \mathbf{V}_{k|k-1}) \) for the state at time \( k \).

\subsection{Kalman Update}

In the influence diagram formulation, the measurement update is performed by introducing measurement evidence and adjusting the diagram’s structure to incorporate this new information. This involves three key steps:

\begin{enumerate}
    \item \textbf{Augmenting the System:}
    \\
    The state and measurement nodes are stacked to form an augmented mean vector and block matrix:
\begin{align}
    \mathbf{u}_{\text{aug}} &=
    \begin{bmatrix}
        \mathbf{u}_{k|k-1} \\
        \mathbf{H} \mathbf{u}_{k|k-1}
    \end{bmatrix}, \\
    \mathbf{V}_{\text{aug}} &=
    \begin{bmatrix}
        \mathbf{V}_{k|k-1} \\
        \mathrm{diag}(\mathbf{R})
    \end{bmatrix}, \\
    \mathbf{B}_{\text{aug}} &=
    \begin{bmatrix}
        \mathbf{B}_{k|k-1}  & \mathbf{H}^\top \\
        \mathbf{0} & \mathbf{0}
    \end{bmatrix}
\end{align}

    where $\mathbf{B}_{k|k-1}$, $\mathbf{V}_{k|k-1}$, and $\mathbf{u}_{k|k-1}$ are the influence diagram parameters for the prior state, and $\mathbf{R}$ is the measurement noise covariance.

    \item \textbf{Evidence Entry:}
    \\
    The measurement $\mathbf{z}$ is entered as evidence for the corresponding measurement node(s). This step updates the state by conditioning the joint Gaussian on the observed value. This is performed by the 'evidence' operation
    \begin{equation}
        \texttt{evidence}\left(\mathbf{u}_{\text{aug}}, B_{\text{aug}}, V_{\text{aug}}, \mathbf{z}, n_0, n_1, n_2, \Delta \mathbf{u}\right)
    \end{equation}
    where $n_0, n_1, n_2$ index the state, measurement, and successor node dimensions, respectively.

    \item \textbf{Arc Reversal:}
    \\
    Evidence entry is implemented via a series of \textit{arc reversals}, which invert the direction of dependencies as needed to ensure the state nodes are conditioned on the measurement. At each step, the reversal operation updates the coefficients and variances according to:
    \begin{align}
        B'_{ji} &= \frac{1}{B_{ij}}, \\
        V'_j &= V_j + \frac{V_i}{B_{ij}^2}
    \end{align}
    for the reversed arc from node $i$ to $j$
\end{enumerate}

After this, the updated mean vector $\mathbf{u}_{k|k}$,  $\mathbf{B}_{k|k}$, and $\mathbf{V}_{k|k}$ reflect the information provided by the measurement $\mathbf{z}$. 

\subsection{Data Association}

For classical JPDAF, data association involves computing association likelihoods for each measurement-to-track pair using the predicted state mean ${\mathbf{x}}_{k|k-1}$ and covariance $\mathbf{P}_{k|k-1}$. Association probabilities are evaluated using the Kalman likelihood: $P(\mathbf{z}_k^{(j)} | {\mathbf{x}}_{k|k-1}) = \mathcal{N}(\mathbf{z}_k^{(j)}; \mathbf{H}_k {\mathbf{x}}_{k|k-1}, \mathbf{S}_k),$
where $\mathbf{S}_k = \mathbf{H}_k \mathbf{P}_{k|k-1} \mathbf{H}_k^\top + \mathbf{R}_k$ is the innovation covariance.

After the measurement update, the $\mathbf{B}_{k|k-1}$ and $\mathbf{V}_{k|k-1}$ are converted back to $\mathbf{P}_{k|k-1}$ to calculate the innovation covariance. But, for gating, the covariance Mahalanobis distance is 
$    d^2 = (\mathbf{z} - \mathbf{z}_{\text{pred}})^\top\, \mathbf{S}^{-1}\, (\mathbf{z} - \mathbf{z}_{\text{pred}})
$.
Using the influence diagram covariance matrix inversion formula \cite{shachter_gaussian_1989}, the  Mahalanobis distance is 
\[
d^2 = (\mathbf{z} - \mathbf{z}_{\text{pred}})^\top
\left(\mathbf{I} - \mathbf{B}\right)^{\top}
\begin{bmatrix}
    \frac{1}{v_1} &        & 0      \\
                  & \ddots &        \\
    0             &        & \frac{1}{v_k}
\end{bmatrix}
\left(\mathbf{I} - \mathbf{B}\right)
(\mathbf{z} - \mathbf{z}_{\text{pred}})
\]
where $\mathbf{V} = [v_j]$
This is calculated without any matrix inversion, making it numerical stable and robust to any ill-conditioning with $S$. 
Overall, the probabilistic association framework (e.g., calculation of association weights, $\beta_k^{(i,j)}$) remains unchanged. The influence diagram structure provides additional numerical stability and flexibility for handling colored noise, model extension, or node removal. 

\section{Simulation Study}
In this section, we compare the performance of ID-JPDAF with JPDAF across several simulation scenarios. We implemented the StoneSoup \cite{hiscocks2023stone} framework to have flexibility to model clutter, detection probability, gating, etc. First, we compare the performance of these algorithms in classical tracking using SIAP metrics \cite{votruba2001single}. Then, we compare the RMSE across a multi-target tracking example in the presence of colored noise. To further investigate filter performance, we tested filter mismatch in the presence of colored noise. Lastly, we test filter performance against varying values of $\sigma$.

\subsection{Tracking Performance under White Noise}

To verify the equivalence between the JPDAF and ID-JPDAF, we simulated a two-dimensional multi-target tracking scenario under Gaussian white noise. Two targets move according to a nearly constant velocity model (NCVM),
\begin{equation}
    \mathbf{x}_{k+1} = \mathbf{F} \mathbf{x}_k + \mathbf{w}_k,
\end{equation}
where $\mathbf{x}_k = [x_{1,k},~ x_{2,k},~ \dot{x}_{1,k},~ \dot{x}_{2,k}]^\top$ denotes the 2-D position and velocity, and $\mathbf{w}_k \sim \mathcal{N}(\mathbf{0},\, \mathbf{Q})$ is an i.i.d.\ process noise sequence. The state transition matrix and process noise covariance are given by
\begin{equation}
    \mathbf{F} = 
    \begin{bmatrix}
        1 & 0 & T & 0 \\
        0 & 1 & 0 & T \\
        0 & 0 & 1 & 0 \\
        0 & 0 & 0 & 1
    \end{bmatrix},
    \quad
    \mathbf{Q} = \sigma_u^2
    \begin{bmatrix}
        \frac{T^3}{3} & 0 & \frac{T^2}{2} & 0 \\
        0 & \frac{T^3}{3} & 0 & \frac{T^2}{2} \\
        \frac{T^2}{2} & 0 & T & 0 \\
        0 & \frac{T^2}{2} & 0 & T
    \end{bmatrix}
\end{equation}
with $T = 1\,\mathrm{s}$ and $\sigma_u^2 = 0.01\,\mathrm{m}^2/\mathrm{s}^4$. The initial state is known, and the prior covariance is set as $\mathbf{P}_0 = 0.01\,\mathbf{I}_4$.

At each time step, each target produces a noisy position measurement,
\begin{equation}
    \mathbf{z}_k = \mathbf{H} \mathbf{x}_k + \mathbf{v}_k,
\end{equation}
where $\mathbf{H} = \begin{bmatrix} 1 & 0 & 0 & 0 \\ 0 & 1 & 0 & 0 \end{bmatrix}$ and $\mathbf{v}_k \sim \mathcal{N}(\mathbf{0},\, \mathbf{R})$ is an i.i.d.\ measurement noise with $\mathbf{R} = \sigma_v^2 \mathbf{I}_2$ and $\sigma_v = 10\,\mathrm{m}$. Each scenario lasts $500$ time steps and is contained in the region of interest (ROI) $[-1000, 1000]\,\mathrm{m} \times [-1000, 1000]\,\mathrm{m}$. Detection probability is $p_\mathrm{d} = 0.5$ and survival probability is $p_\mathrm{s} = 0.995$. Clutter measurements are Poisson-distributed with mean $\mu_c = 5$ and uniformly distributed over the ROI.

Both the classical JPDAF and the ID-JPDAF were run using this simulation setup. Over $1000$ Monte Carlo trials, the root mean squared error (RMSE) and estimated covariances produced by both methods agreed to within $10^{-6}$ at every time step, confirming the numerical equivalence of the two implementations under white noise.

\subsection{Tracking Performance under Colored Noise}
We evaluated filter performance in a challenging multi-target tracking scenario characterized by temporally correlated (AR(1)) measurement noise. The scenario involved two targets moving in a two-dimensional plane, each governed by an NCVM model. The state vector for each target was defined as $\mathbf{x}_k = [x_k, \dot{x}_k, y_k, \dot{y}_k, n_{x,k}, n_{y,k}]^\top$, where $(x_k, y_k)$ and $(\dot{x}_k, \dot{y}_k)$ denote position and velocity, and $(n_{x,k}, n_{y,k})$ represent auxiliary states for the colored noise process.

$\mathbf{w}_k$ is zero-mean Gaussian process noise with covariance $\mathbf{Q} = \mathrm{diag}(q_p, q_v, q_p, q_v, q_n, q_n)$. Here, $q_p = 100$, $q_v = 10$, $q_n = 25$.

The AR(1) process for the colored noise terms followed
\begin{equation}
    n_{x, k+1} = \rho n_{x, k} + \xi_{x,k}, \qquad n_{y, k+1} = \rho n_{y, k} + \xi_{y,k},
\end{equation}
with noise $\xi_{x,k}, \xi_{y,k} \sim \mathcal{N}(0, \sigma^2)$, where $\sigma = 3.0$m.

Measurements were generated at each time step via an \emph{augmented measurement model}:
\begin{equation}
    \mathbf{z}_k = 
    \begin{bmatrix}
        x_k + n_{x,k} \\
        y_k + n_{y,k}
    \end{bmatrix}
\end{equation}
which corresponds to the standard position measurements corrupted by AR(1) noise.

Simulation experiments were conducted for $2000$ time steps of length $\tau$ = $0.05s$, repeated over $200$ random seeds. The initial target states were set to $[0, 1, 0, 1, 0, 0]^\top$ and $[20, -1, 10, -1, 0, 0]^\top$ for the two targets. Filters were initialized with a prior covariance $\mathbf{P}_0 = \mathrm{diag}(100, 10, 100, 10, 25, 25)$. 

\begin{figure}[!t]
    \centering
    \includegraphics[width=0.48\textwidth]{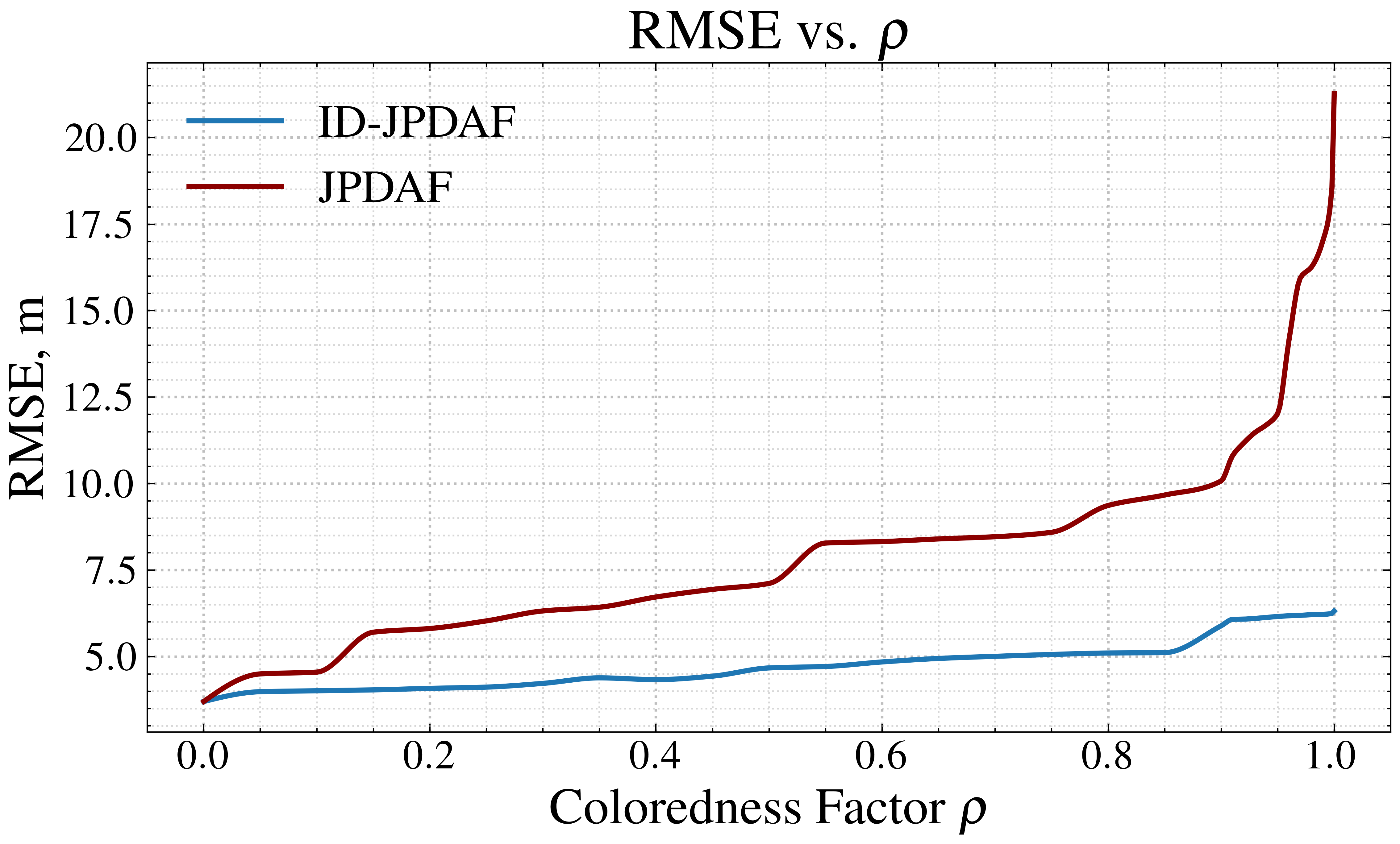}
    \caption{
        This figure plots the value of RMSE with the aforementioned simulation setup, where the blue line represents ID-JPDAF and the red line represents JPDAF
    }
    \label{fig:rmse_vs_rho}
\end{figure}

Performance was evaluated using the RMSE between estimated and ground truth trajectories. Results were averaged over all Monte Carlo trials. First, we calculated the RMSE over different values of $\rho$ to investigate the effect of colored noise on filter performance. Figure \ref{fig:rmse_vs_rho} shows the evolution of mean position RMSE as a function of $\rho$. For low‐correlation ($\rho < 0.2$), both filters perform comparably when the noise is nearly white with RMSE $ \simeq4$ m. As $\rho$ increases to 0.2, JPDAF’s error rises slightly faster ($ \simeq 5.5$ m) than ID‑JPDAF ($\simeq 4.2m$), indicating that even modest temporal correlations begin to degrade association quality. For moderate correlation ($0.2 < \rho < 0.8$). JPDAF’s RMSE climbs steadily from $\simeq 6$ m to $\simeq 9$ m, while ID‑JPDAF remains below $5$ m throughout. 
For high-correlation ($\rho > 0.8$), JPDAF’s error accelerates sharply—reaching over $12$ m at $\rho$ = $0.9$ and increases asymptotically as $\rho \to 1$. However, ID‑JPDAF only increases modestly to about $6.3$ m.  Overall, ID‑JPDAF yields a considerable reduction in RMSE for even moderate correlation, and yields a significant reduction in RMSE for high correlation.

\begin{figure}[!t]
    \centering
    \includegraphics[width=0.48\textwidth]{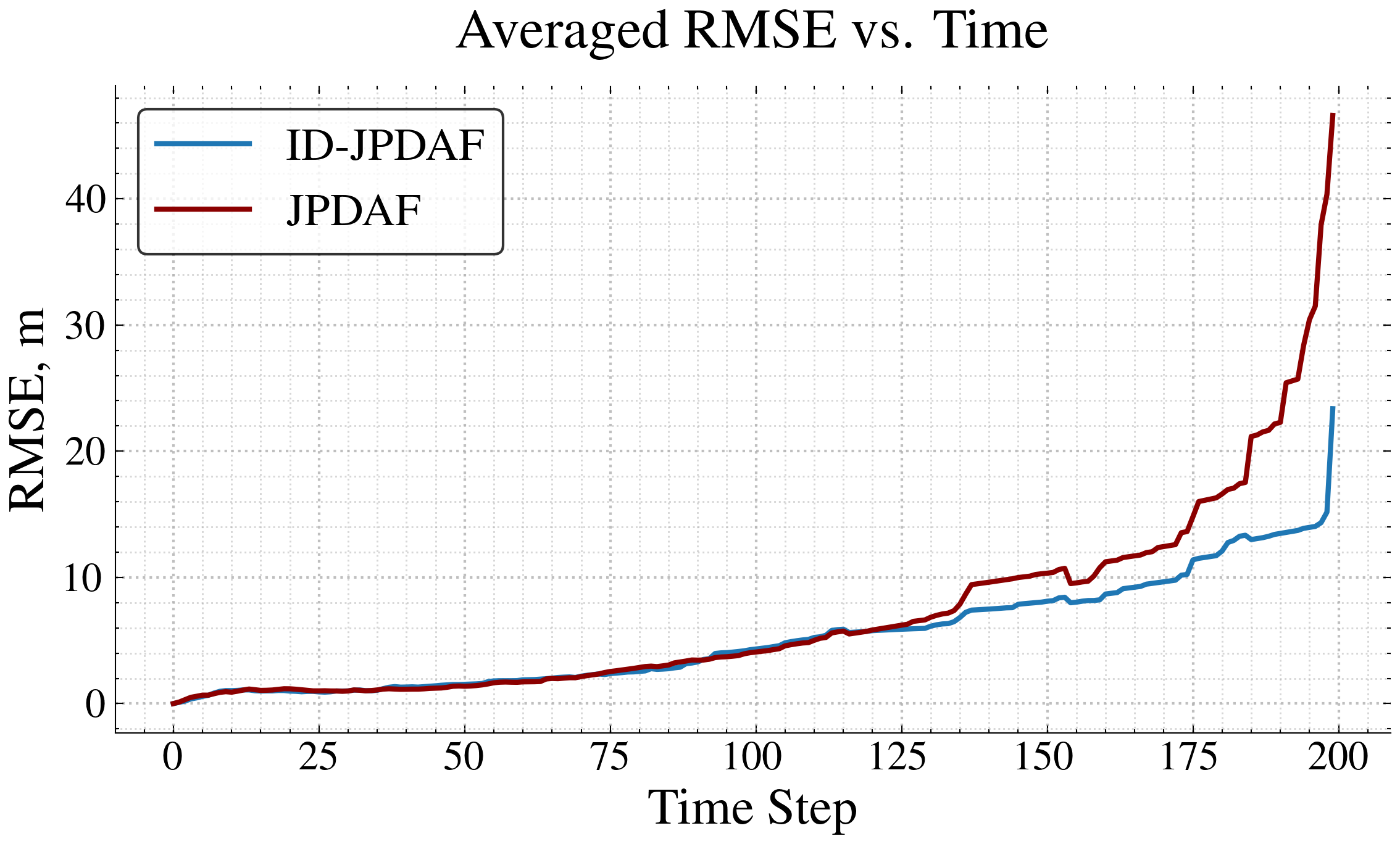}
    \caption{RMSE with the filter mismatch where true values are underestimated.}
    \label{fig:filter_mismatch}
\end{figure}

Then, we investigate two filter mismatch scenarios. The first case underestimates the true values of $\rho$ and $\sigma$, while the second case overestimates 4these parameters. For the first case, $\rho_{true}$ = $0.95$ and $\sigma_{true} = 1.5$ while the filter assumes $\rho_{filter}$ = $0.9$ and $\sigma_{filter} = 0.5$. Figure \ref{fig:filter_mismatch} plots the time‑evolution of the mean position RMSE for both the filters. For the first 120 frames, both filters track almost identically, with RMSE climbing gently from $0$ up to $\simeq 6$ m. Beyond Step 120, JPDAF’s error grows noticeably faster, hitting $\simeq 10$ m around step 140 and then accelerating sharply after step 170 to exceed $45$ m by step 200. In contrast, ID‑JPDAF remains far more stable, rising only to $\simeq 12$ m at step 200. 4

\begin{figure}[!t]
    \centering
    \includegraphics[width=0.48\textwidth]{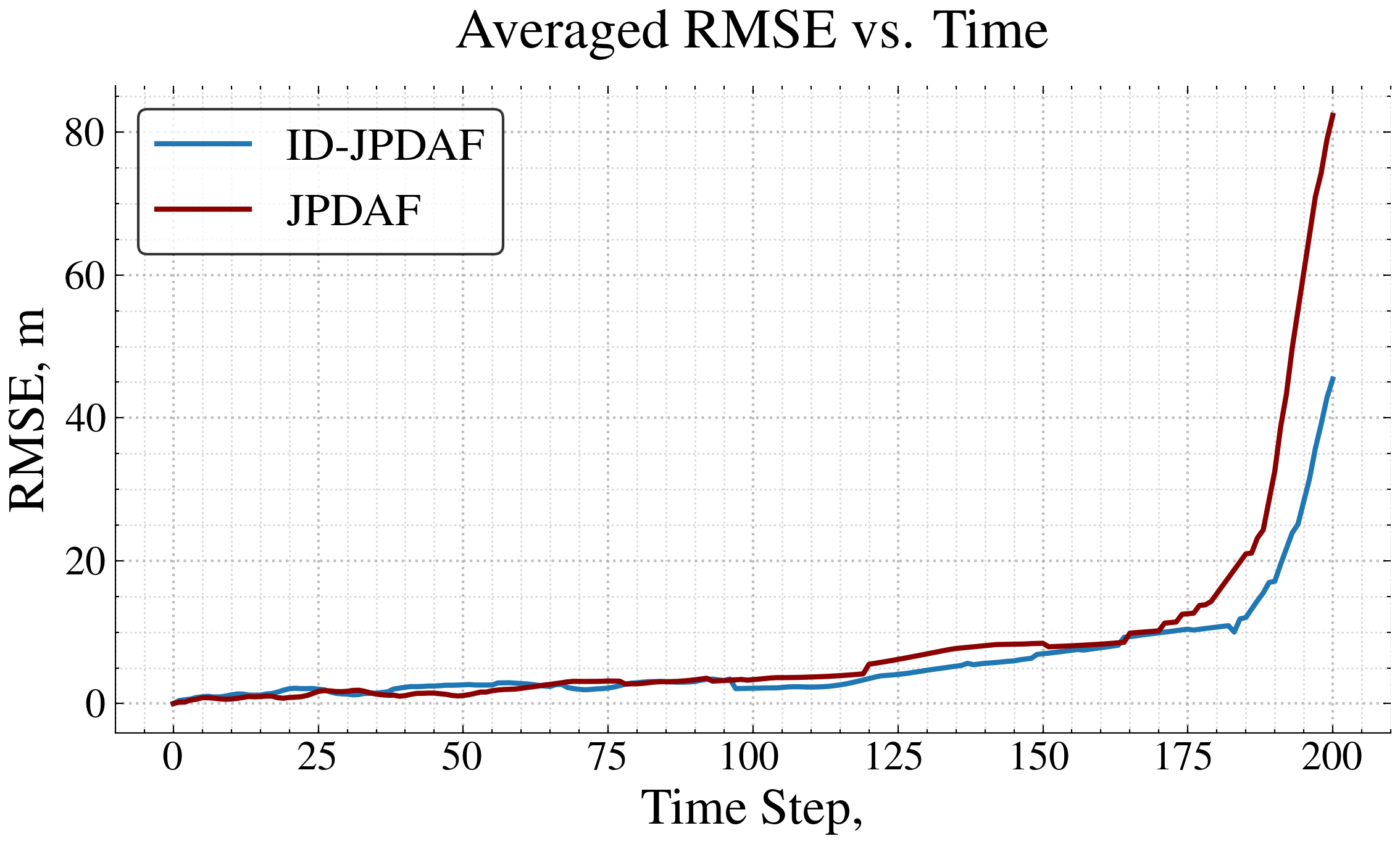}
    \caption{
        RMSE with the filter mismatch where true values are overestimated.
    }
    \label{fig:filter_mismatch_2}
\end{figure}

For the second case, $\rho_{true}$ = $0.9$ and $\sigma_{true} = 0.5$ while the filter assumes $\rho_{filter} = 0.95$ and $\sigma_{filter} = 1.5$ Figure \ref{fig:filter_mismatch_2} plots the time‐evolution of mean position RMSE for both filters. Both filters remain tightly coupled for the first $175$ frames, but beyond this JPDAF’s error accelerates much more rapidly, reaching $80$ by $t=200$ and blowing up to $80$ m. On the other hand, ID‑JPDAF degrades far more gracefully, only climbing to $45$ m at
$t=200$. 

\begin{figure}[!t]
    \centering
    \includegraphics[width=0.48\textwidth]{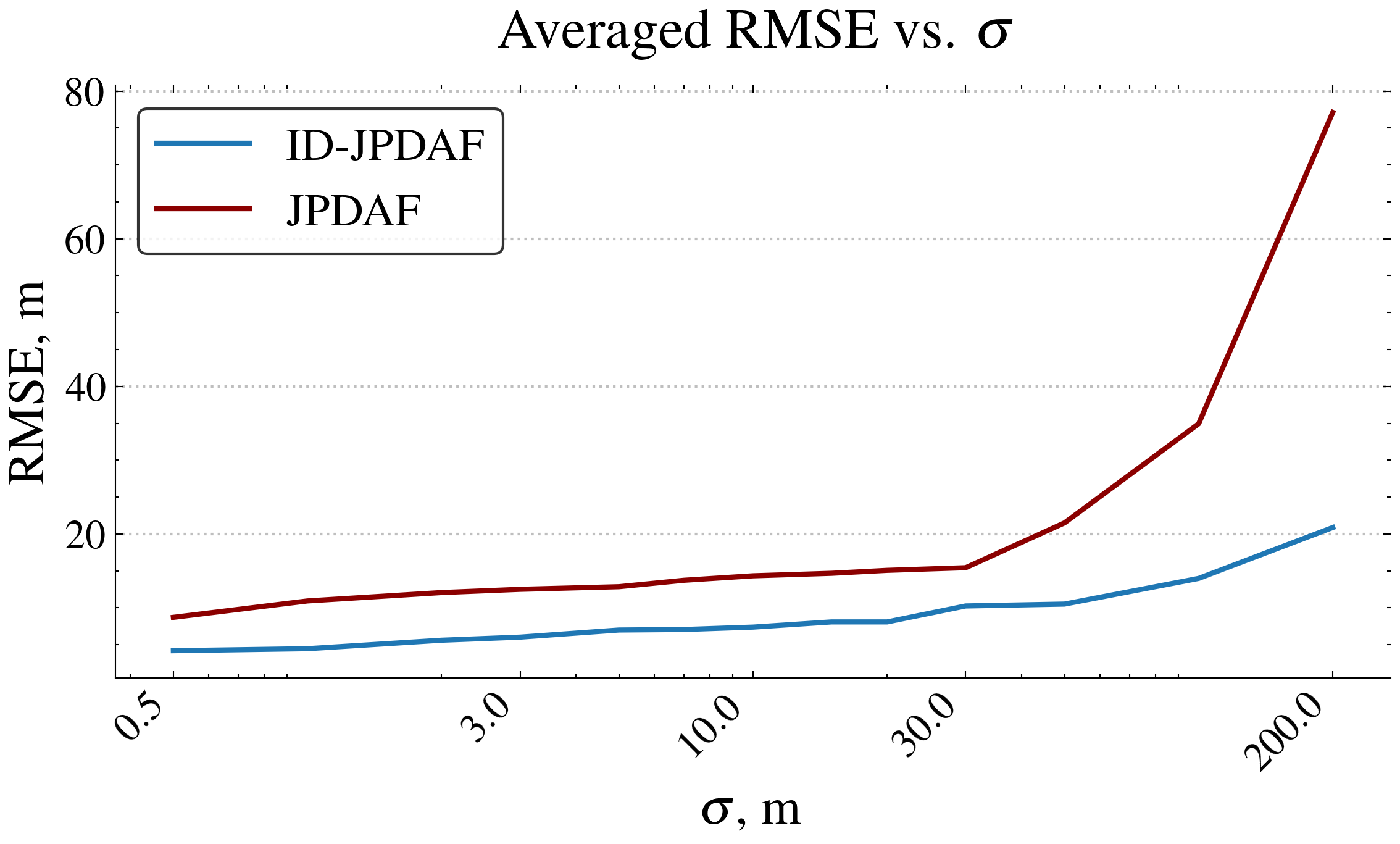}
    \caption{
        This figure plots RMSE across different values of $\sigma$. $\rho = 0.80$.
    }
    \label{fig:RMSE_vs_Sigma}
\end{figure}

Lastly, we investigate filter performance for different values of $\sigma$ where $\rho = 0.80$. Figure \ref{fig:RMSE_vs_Sigma} shows RMSE of the filters as $\sigma$ increases from $0.5$ m to $200$ m. Both methods start at roughly $5$ m RMSE when $\sigma = 0.5$ m, but as $\sigma$ grows, JPDAF’s error climbs faster at $\sigma =30$ m. Whereas ID‑JPDAF increases more gently, staying below $20$ m throughout the entire process. 

\section{Conclusions}
In this paper, we have identified and analyzed two fundamental limitations of classical Kalman filter‐based JPDAF in practical multi‐target tracking: sensitivity to colored measurement noise and numerical instability under ill‐conditioned covariance updates. To address these issues, we introduced an Influence Diagram–based JPDAF (ID‐JPDAF) that embeds the Gaussian influence‐diagram formalism into the Kalman prediction and update steps, yielding a numerically robust filtering procedure without explicit covariance inversions.

Through extensive StoneSoup–based simulations, we demonstrated that under moderate to strong noise correlation ($\rho>0.2$), ID‐JPDAF achieves a 30–70\% reduction in RMSE compared to standard JPDAF, and remains stable even as $\rho\to1$. This improvement happens due to the conversion of covariance matrix $\mathbf{P}$ to regression coefficients $\mathbf{B}$ and conditional variance $\mathbf{V}$. Since this better encodes the conditional dependence across operations, ID-JPDAF performs better in the presence of colored noise. Further, since colored measurement noise can cause numerical instability due to ill-conditionedness \cite{maybeck1982stochastic}, the lack of inversion in ID-JPDAF results in lower RMSE. With Bayesian network models, this conditional breakdown aids filter performance \cite{ghosh2008probabilistic, koller2009probabilistic}. 

In scenarios involving filter-model mismatch, standard JPDAF exhibits significant divergence, whereas ID-JPDAF demonstrates graceful degradation, limiting peak RMSE to less than half that of JPDAF. This performance advantage holds in both underestimation and overestimation of model parameters. Furthermore, across a wide range of measurement noise variances ($\sigma = 0.5$ m $\to$ $200$ m), ID-JPDAF consistently outperforms JPDAF, exhibiting greater robustness and stability.

These results confirm that influence‐diagram inference offers a powerful means to enhance multi‐target tracking robustness in realistic radar environments. Future work will explore extensions to nonlinear dynamics, adaptive tuning of influence‐diagram parameters, and performing analysis on the scalability of ID-JPDAF.

\bibliographystyle{IEEEtran}
\bibliography{references}

\begin{thebibliography}{10}
\providecommand{\url}[1]{#1}
\csname url@samestyle\endcsname
\providecommand{\newblock}{\relax}
\providecommand{\bibinfo}[2]{#2}
\providecommand{\BIBentrySTDinterwordspacing}{\spaceskip=0pt\relax}
\providecommand{\BIBentryALTinterwordstretchfactor}{4}
\providecommand{\BIBentryALTinterwordspacing}{\spaceskip=\fontdimen2\font plus
\BIBentryALTinterwordstretchfactor\fontdimen3\font minus
  \fontdimen4\font\relax}
\providecommand{\BIBforeignlanguage}[2]{{%
\expandafter\ifx\csname l@#1\endcsname\relax
\typeout{** WARNING: IEEEtran.bst: No hyphenation pattern has been}%
\typeout{** loaded for the language `#1'. Using the pattern for}%
\typeout{** the default language instead.}%
\else
\language=\csname l@#1\endcsname
\fi
#2}}
\providecommand{\BIBdecl}{\relax}
\BIBdecl

\bibitem{kalman1960}
R.~E. Kalman, ``A new approach to linear filtering and prediction problems,''
  \emph{Journal of Basic Engineering}, vol.~82, no.~1, pp. 35--45, 1960.

\bibitem{maybeck1982stochastic}
P.~S. Maybeck, \emph{Stochastic models, estimation, and control}.\hskip 1em
  plus 0.5em minus 0.4em\relax Academic press, 1982, vol.~3.

\bibitem{blackman1986multiple}
S.~S. Blackman, ``Multiple-target tracking with radar applications,''
  \emph{Dedham}, 1986.

\bibitem{chang1984joint}
K.-C. Chang and Y.~Bar-Shalom, ``Joint probabilistic data association for
  multitarget tracking with possibly unresolved measurements and maneuvers,''
  \emph{IEEE Transactions on Automatic control}, vol.~29, no.~7, pp. 585--594,
  1984.

\bibitem{reid2003algorithm}
D.~Reid, ``An algorithm for tracking multiple targets,'' \emph{IEEE
  transactions on Automatic Control}, vol.~24, no.~6, pp. 843--854, 2003.

\bibitem{kenley1986influence}
C.~R. Kenley, \emph{INFLUENCE DIAGRAM MODELS WITH CONTINUOUS VARIABLES
  (DECISION ANALYSIS, MULTIVARIATE NORMAL DISTRIBUTION,
  LINEAR-QUADRATIC-GAUSSIAN, KALMAN FILTERS, COVARIANCE ASSESSMENT)}.\hskip 1em
  plus 0.5em minus 0.4em\relax Stanford University, 1986.

\bibitem{bryson_chapter_1975}
A.~E. Bryson and Y.-C. Ho, ``\BIBforeignlanguage{English}{Chapter 12: {Optimal}
  filtering and prediction},'' in \emph{\BIBforeignlanguage{English}{Applied
  optimal control: optimization, estimation, and control}}, revised
  printing~ed.\hskip 1em plus 0.5em minus 0.4em\relax Washington: Hemisphere,
  1975, pp. 348--389.

\bibitem{barshalom1993}
Y.~Bar-Shalom and X.~R. Li, \emph{Estimation and Tracking: Principles,
  Techniques, and Software}.\hskip 1em plus 0.5em minus 0.4em\relax Artech
  House, 1993.

\bibitem{michael_kovacich_application_1990}
\BIBentryALTinterwordspacing
{Michael Kovacich}, ``Application of {Bayesian} networks to multitarget
  tracking,'' vol. 1305, Oct. 1990, p. 348. [Online]. Available:
  \url{https://doi.org/10.1117/12.2321776}
\BIBentrySTDinterwordspacing

\bibitem{howard2005influence}
R.~A. Howard and J.~E. Matheson, ``Influence diagrams,'' \emph{Decision
  Analysis}, vol.~2, no.~3, pp. 127--143, 2005.

\bibitem{shachter_gaussian_1989}
\BIBentryALTinterwordspacing
R.~D. Shachter and C.~R. Kenley, ``Gaussian {Influence} {Diagrams},''
  \emph{Management Science}, vol.~35, no.~5, pp. 527--550, 1989. [Online].
  Available: \url{http://www.jstor.org/stable/2632102}
\BIBentrySTDinterwordspacing

\bibitem{hiscocks2023stone}
S.~Hiscocks, J.~Barr, N.~Perree, J.~Wright, H.~Pritchett, O.~Rosoman,
  M.~Harris, R.~Gorman, S.~Pike, P.~Carniglia \emph{et~al.}, ``Stone soup: No
  longer just an appetiser,'' in \emph{2023 26th International Conference on
  Information Fusion (FUSION)}.\hskip 1em plus 0.5em minus 0.4em\relax IEEE,
  2023, pp. 1--8.

\bibitem{votruba2001single}
P.~Votruba, R.~Nisley, R.~Rothrock, and B.~Zombro, ``Single integrated air
  picture (siap) metrics implementation,'' Tech. Rep., 2001.

\bibitem{ghosh2008probabilistic}
J.~K. Ghosh, ``Probabilistic networks and expert systems: Exact computational
  methods for bayesian networks by robert g. cowell, a. philip dawid, steffen
  l. lauritzen, david j. spiegelhalter,'' 2008.

\bibitem{koller2009probabilistic}
D.~Koller and N.~Friedman, \emph{Probabilistic graphical models: principles and
  techniques}.\hskip 1em plus 0.5em minus 0.4em\relax MIT press, 2009.

\end{thebibliography}
\end{document}